\documentclass[12pt]{iopart}

\begin{document}

\title[Group-Theoretic Matching of the Length and Equality Principles in Geometry]
{Group-Theoretic Matching of the Length and Equality Principles in Geometry}

\author{Serhii Samokhvalov, Olena Balakireva}

\address{Dniprovsk State Technical University, Kamianske, UKRAINE}
\ead{serg\_samokhval@ukr.net}
\vspace{10pt}

\begin{abstract}
Deformed generalized gauge groups, which were created from physical considerations and made it possible to clarify some long-standing problems in physics, such as the problem of motion and the problem of the energy of the gravitational field, simultaneously carried out the implementation of the Klein’s Erlangen program for spaces with variable curvature, and for Riemannian spaces even in two different ways.
In this paper, this issue is considered from a geometric point of view and these two methods of group-theoretical description of Riemannian spaces are reconciled.
The paper deals with the canonical deformed group of diffeomorphisms with a given length scale which describes the motion of unit scales in a Riemannian space. This allows one to measure the lengths of arbitrary curves implementing the length principle which was laid by B. Riemann at the foundation of geometry. We present a method of univocal extension of this group to a group which contains gauge rotations of vectors (the group of parallel transports) whose transformations leave unchanged the lengths of vectors and the corners between them, implementing for Riemannian spaces the Klein's principle of equality and matching both principles of the foundations of geometry, thus overcoming the Riemann-Klein antagonism.

\end{abstract}

%
\vspace{2pc}

\section*{Introduction}

Creation of Lobachevskii geometry and Gauss’ works on the theory of surfaces stimulated the intensive development of geometric theories in the 19-th century. As a result, the works of B. Riemann \cite{1} and F. Klein \cite{2} were appeared. These works laid different principles at the foundation of the geometry: the \textit{length principle} which requires the possibility to measure the lengths of arbitrary lines no matter how they are situated, and \textit{equality principle} which is established  by coincidence of figures in the space by means of transformations belonging to a group of transformations of the space – the principal group of the geometry under consideration (according to F. Klein). According to E. Cartan \cite{3}, there is an \textit{antagonism} between these two principles owning to the absence of any homogeneity in an arbitrary curved Riemannian space.

Earlier, attempts were made to overcome this antagonism by means of refusing from group structure of used transformations. E. Cartan suggested to consider a curved space as a \textit{anholonomic space} with the same principal (fundamental according to the terminology used by Cartan) group as the corresponding flat space \cite{3}. R. Sulanke and P. Wintgen \cite{4} applied the \textit{category theory} for description of curved spaces, and L. V. Sabinin \cite{5} used \textit{quasigroups} replying on the thesis that nonassociativity is an algebraic equivalent of the geometric notion of curvature.

It has been clarified later that so-called deformed generalized gauge groups introduced from physical arguments \cite{6} can describe, and in two different ways, geometric structure of variable curvature \cite{7}, in particular, affinely connected spaces \cite{8} and Riemannian spaces \cite{9}.

In the first case, the translation generators of a group coincide with covariant derivatives in a curved affinely connected space. The curvature and the torsion of the space are determined by the structure functions of the group, which are antisymmetric parts of the coefficients of the second order terms in the expansion of the multiplication law of the group with respect to parameters \cite{8}.

In the second case, it has become possible to use for torsion-free affinely connected spaces, in particular, for Riemannian spaces, a more narrow group, namely, the deformed group of diffeomorphisms. Acting on a manifold, it generates an action on its tangent bundle and thus defines a rule of parallel translations of vectors. But, in this case, for description of the curvature tensor, one needs to know the third terms of the expansion of the multiplication law of the group with respect to parameters \cite{9}.

In this paper, we suggest a method which makes it possible to match and unite these two approaches to group-theoretical description of Riemannian spaces and thus to overcome the Riemann-Klein antagonism.

The method we use eliminates the indefiniteness which took place in previous papers \cite{8},\cite{9} with respect to finite translations of vectors given in curved spaces by action of deformed generalized gauge groups. This has been achieved by means of a generalization to the case under consideration of the notion of canonical Lie groups, which establishes one-to-one correspondence between the groups and the geometric structures defined by them. We present criteria of canonicity of generalized gauge groups under consideration.

We show that a canonical deformed group of diffeomorphisms with given length scale – the group of Riemannian translations $RT$ – makes it possible to measure the lengths of geodesics by means of translation of the unit scale along itself, which realizes the \textit{length principle} (Riemann’s approach). What is more, the group $RT$ contains in itself all characteristics of Riemannian space, in particular, the \textit{curvature as a characteristic of deviation of geodesics}, and defines by its action on a manifold a \textit{Riemannian structure} (by means of geometric objects of a Riemannian geometry such as vector and orthonormal frame fields, and, therefore, the metric). But the group $RT$ does not preserve invariants on a set, and, therefore, cannot be considered as the principal group of a Riemannian space.

However, as is shown in this paper, a given Riemannian structure on a manifold (a metric) can be defined by different groups $RT$ which are related to each other by nonlocal automorphisms of special form. The set of such automorphisms is a group which is, in the case the requirement of its canonicity is imposed, the group $DP$ of parallel translations of vectors in a Riemannian space. The group $DP$, acting in the tangent bundle of a manifold, preserves the lengths of vectors and the angles between vectors. Thus, the group of parallel transports $DP$ can be considered as the \textit{principal group of a Riemannian space} because it does not change its characteristics and enables ones to establish the \textit{equality} (congruence) of vectors and their mutual disposition even in the case they are given in different points (Klein’s approach). Here \textit{the curvature plays the role of noncommutativity of covariant derivatives} – the generators of translations of $DP$ (or of the rotations of a vector when it is displaced along a closed contour). Acting in the bundle of orthonormal frames on a manifold (a principal bundle), the group $DP$ describes motions of frames, which is an interpretation of E. Cartan’s \textit{method of moving frames} \cite{10}, though the structure equation of s Riemannian space in this case is not postulated as it is within Cartan’s approach, but is a consequence of the existence of a group $DP$.

The group of parallel transports $DP$, which realizes for Riemannian space the equality principle, has as a subgroup the group of Riemannian translations $RT$, which realizes the length principle and, therefore, unites the two approaches laid at the foundation of geometry by B. Riemann and F. Klein, thus \textit{overcoming the Riemann-Klein antagonism}.

In the paper, we do not consider global-topological questions, all relations are obtained within a coordinate domain of manifold, though we use more adopted in contemporary geometry approach which avoids the use of coordinates. In addition, by groups we mean the corresponding local groups, and some terminology and notations concerning Lie groups follows \cite{11}.

\section{Deformed gauge group of translations}

\textit{The gauge group of translations} $T_X^g$ of an affine space $(X,\widetilde{T})$ whose elements are vector functions $\widetilde{t}_x$ with values in the vector space $\widetilde{T}$ which depend smoothly on $x\subset X$ and satisfy the condition $|1+\partial_x \widetilde{t}_x|\neq 0$, where $\partial_x := \partial/\partial x$, has the following multiplication law:

\begin{equation}\label{1}
(\widetilde{t} \times \widetilde{t'})_x=\widetilde{t}_x+\widetilde{t'}_{x'}, \quad x'=x+\widetilde{t}_x, \quad \widetilde{t}_x,\,\widetilde{t'}_{x'}\subset T_X^g.
\end{equation}
Thus, $T_X^g$ is group of diffeomorphisms of the space $X$ in additive parameterization with parameters $\widetilde{t}_x$.

We see that the multiplication law in the group $T_X^g$ gives a rule for the addition of vectors specified at different points, which is a characteristic feature of generalized gauge groups, which is manifested in their important role in describing geometric structures in which parallel transport takes place.

The subset of elements of the group $T_X^g$, which is parameterized by vectors independent of $x$, obviously forms the group $\widetilde{T}$.

\textit{Deformed gauge translation group} $T_X^{gH}$ (deformed group of diffeomorphisms) is by definition isomorphic to the group $T_X^g$ and is parametrized by vector functions $t_x = H_x (\widetilde{t}_x)$ taking values in a vector space $T$ isomorphic to $\widetilde{T}$, where $H_x : \widetilde{T} \rightarrow T$ is an invertible smoothly depending on $x, \: \forall x \in X$ a \textit{deformation map}, $K_x := H_x^{-1}$. The mapping $H_x$ is generally nonlinear and, therefore, is not an isomorphism between the vector spaces $\widetilde{T}$ and $T$.

The multiplication law and the action of the group $T_X^{gH}$ on $X$ are determined by its isomorphism to the group $T_X^g$ \cite{9}:

\begin{equation} \label{2}
\eqalign{ {(t\times t')_{x}=\varphi_{x}(t_{x},t'_{x'}):=H_{x}(K_{x}(t_{x})+K_{x'}(t'_{x'})),  } \cr
{x'=f_{x}(t_{x}):=x+K_{x}(t_{x}). } } \end{equation}

The associativity law of multiplication and the composition law of the action of the group $T^{gH}_{X}$ on $X$ are expressed by the following equations for the mappings $\varphi_{x}$ and $f_{x}$:

\begin{equation}\label{3}
\varphi_{x}(\varphi_{x}(t,t'),t'')=\varphi_{x}(t,\varphi_{x'}(t',t'')),\quad f_{x'}(t')=f_{x}(\varphi_{x}(t,t')),
\end{equation}
which are fulfilled due to the deformation method of constructing the group $T^{gH}_{X}$ and whose solution is equalities (\ref{2}) with an arbitrary deformation mapping $H_{x}$ and constant parameters. The subset of elements of the deformed group of diffeomorphisms $T^{gH}_{X}$, parametrized by vectors independent of $x$, for an arbitrary deformation, in contrast to the undeformed case, is not closed with respect to the multiplication operation, and therefore no longer forms a group.

	A consequence of equations (\ref{3}) are equations that generalize the Lie equations of the theory of Lie groups for the group $T^{gH}_{X}$:

\begin{equation}\label{4}
\mu_{x}(t)\cdot\partial\varphi_{x}(t,t')=\mu_{x}(\varphi_{x}(t,t'))+e_{x}\varphi_{x}(t,t'),
\end{equation}
\begin{equation}\label{5}
\lambda_{x'}(t')\cdot\partial'\varphi_{x}(t,t')=\lambda_{x}(\varphi_{x}(t,t')),
\end{equation}
\begin{equation}\label{6}
\mu_{x}(t)\cdot\partial f_{x}(t)=e_{x}f_{x}(t),\quad\lambda_{x}(t)\cdot\partial f_{x}(t)=k_{x'},
\end{equation}
where $\mu_{x}(t'):=\partial\varphi_{x}(0,t')$, $\lambda_{x}(t):=\partial'\varphi_{x}(t,0)$ are \textit{auxiliary maps}, $e_{x}:=k_{x}\cdot\partial_{x}$ are the \textit{generators} of the action of the group $T^{gH}_{X}$ on $X$, $k_{x}:=\partial K_{x}(0)$ are the \textit{deformation coefficients} (here it is assumed $\partial:=\partial/\partial t$, $\partial':=\partial/\partial t'$).

The integrability conditions for equations (\ref{4}) - (\ref{6}) are the equations that generalize for the group $T^{gH}_{X}$ the \textit{Maurer – Cartan equations} of the theory of Lie groups:

\begin{equation}\label{7}
e_{x}\mu_{x}(t)\langle l,l'\rangle+\partial\mu_{x}(t)\langle\mu_{x}(t)\langle l'\rangle,l\rangle
-(l\leftrightarrow l')=\mu_{x}(t)\langle C_{x}\langle l,l'\rangle\rangle,
\end{equation}
\begin{equation}\label{8}
\partial\lambda_{x}(t)\langle\lambda_{x}(t)\langle l\rangle,l'\rangle
-(l\leftrightarrow l')=\lambda_{x}(t)\langle C_{x'}\langle l,l'\rangle\rangle,
\end{equation}
\begin{equation}\label{9}
e_{x}k_{x}\langle l,l'\rangle-(l\leftrightarrow l')=k_{x}\langle C_{x}\langle l,l'\rangle\rangle,
\end{equation}
where $l,l'\in T$ and

\begin{equation}\label{10}
C_{x}\langle l,l'\rangle:=\gamma_{x}\langle l,l'\rangle-(l\leftrightarrow l')
\end{equation}
is the \textit{structure operator} of the group $T^{gH}_{X}$ that depends explicitly on $x$ (in contrast to the case of finite-dimensional Lie groups), which is the antisymmetric part of the mapping $\gamma_{x}:=\partial\partial'\varphi_{x}(0,0)$, which defines the second order of the group $T^{gH}_{X}$  multiplication law and depends explicitly on $x$ due to the explicit the dependence of the mapping $\varphi_{x}$ of the group $T^{gH}_{X}$ on $x$.

The integrability condition for equations (\ref{7}) - (\ref{9}) is the \textit{Jacobi relation} for the structure operator, which in our case takes the form

\begin{equation}\label{11}
e_{x}C_{x}\langle l,l',l''\rangle+C_{x}\langle l,C_{x}\langle l', l''\rangle\rangle
+cycl(l,l',l'')=0.
\end{equation}
	
Differentiation of equation (\ref{7}) with respect to $t$ at zero allows us to reveal an important relation

\begin{equation}\label{12}
e_{x}\gamma_{x}\langle l,l',t\rangle+\gamma_{x}\langle l,\gamma_{x}\langle l',t\rangle\rangle
-(l\leftrightarrow l')=R_{x}\langle t,l,l'\rangle+\gamma_{x}\langle C_{x}\langle l,l'\rangle ,t\rangle
\end{equation}
for the antisymmetric part

\begin{equation}\label{13}
R_{x}\langle t,l',l\rangle:=\rho_{x}\langle l,l',t\rangle-(l\leftrightarrow l')
\end{equation}
of the mapping $\rho_{x}:=\partial\partial'^{2}\varphi_{x}(0,0)$, which partially determines the third order of the multiplication law for the group $T^{gH}_{X}$. The mapping $R_x$ will be called the \textit{curvature operator} of the group $T^{gH}_{X}$.

Definition (\ref{13}) and the symmetry of the mapping $\rho_{x}$ in the last two arguments imply the cyclic identity

\begin{equation}\label{14}
R_{x}\langle l,l',l''\rangle+cycl(l,l',l'')=0,
\end{equation}
which, due to Eq. (\ref{12}), reduces to the Jacobi relation (\ref{11}).

From a geometric point of view, the generators $e_x$ of the group $T^{gH}_{X}$ define a \textit{field of affine frames} on $X$, and its structure operator $C_{x}$, due to Eq. (\ref{9}), is its \textit{anholonomic object}.

The multiplication law for the group $T^{gH}_{X}$, written for the infinitesimal parameters $(t\times\theta)_{x}=
t_{x}+\theta_{x'}+\gamma_{x}\langle t_{x},\theta_{x'}\rangle$, gives a rule for the composition of the vectors $t_x$ and $\theta_{x'}$ given at different points $x$ and $x'$, and therefore and a certain \textit{rule for parallel transport} of the vector $\theta_{x'}$ from point $x'$ to point $x$:

\begin{equation}\label{15}
\theta_{x||}:=\theta_{x'}+\gamma_{x}\langle t_{x},\theta_{x'}\rangle=\theta_{x}+t_{x}\cdot \nabla_{x}\theta_{x},
\end{equation}
with which the composition law of the group $T^{gH}_{X}$ describes the addition of vectors given already at one point: $(t\times\theta)_{x}=t_{x}+\theta_{x||}$. Here $\nabla_{x}=e_{x}+\gamma_{x}$ is the covariant derivative, and in this sense the mapping $\gamma_{x}$ plays the role of an \textit{object of affine connection} in the frame $e_{x}$. Due to Eq. (\ref{10}), there is no torsion in our case.

The curvature operator $R_x$ of the group $T^{gH}_{X}$, by virtue of Eq. (\ref{12}), which in this case plays the role of the \textit{structure equation of the affinely connected space}, is its \textit{curvature tensor}.

Using (\ref{2}), we find

\begin{equation}\label{16}
\gamma_{x}\langle l,l'\rangle=h_{x}\langle\Gamma_{x}\langle k_{x}\langle l\rangle,k_{x}\langle l'\rangle\rangle
+e_{x}k_{x}\langle l,l'\rangle\rangle,
\end{equation}
where $h_{x}:=k^{-1}_{x}$ and $\Gamma_{x}:=k_{x}\circ\tilde{\partial}^{2}H_{x}(0)$, here $\tilde{\partial}:=\partial/\partial\tilde{t}$. Thus, the covariant derivative can be represented as $\nabla_{x}=h_{x}\circ k_{x}\cdot\tilde{\nabla}_{x}\circ k_{x}$, where $\tilde{\nabla}_{x}=\partial_{x}+\Gamma_{x}$ is the covariant derivative in the frame $\partial_{x}$. Consequently, a smooth symmetric bilinear map $\Gamma_{x}$ is an object of torsion-free affine connection in the holonomic frame $\partial_{x}$, and since it is other-wise arbitrary due to the arbitrariness of the deformation mappings $H_x$, the second order of which is given by $\Gamma_{x}$, an arbitrary torsion-free affine connection can be described in this way.

\textbf{Theorem 1}. \textit{Acting on a manifold $X$, the deformed gauge group of translations $T^{gH}_{X}$ defines on it the structure of an affinely connected torsion-free space whose structure equation is a necessary condition for the existence of the group $T^{gH}_{X}$. An arbitrary torsion-free affine connection on $X$ can be defined in this way.}

Obviously, the group $T^{g}_{X}$  defines on $X$ the structure of a flat affine space.

\section{Canonical deformations}

If in the expansion of the deformation map

\[
H_{x}(\tilde{t})=h_{x}\circ(\tilde{t}+\frac{1}{2}\Gamma_{x}\langle\tilde{t},\tilde{t}\rangle+\ldots),
\]
with the help of which the group $T^{gH}_{X}$ is constructed, the first order determines the field of affine frames on $X$, and the second - the affine connection, then all other orders do not affect to the connection defining by the group $T^{gH}_{X}$ and can be arbitrary. To remove this ambiguity, we generalize the notion of canonicity of Lie groups for the groups $T^{gH}_{X}$.

A Lie group $G$, $g\in G$, is \textit{canonical} if any line of the form $g(s)=s\tau$ is its one-parameter subgroup \cite{11}. In this case, the mapping $\varphi(g,g'):=g\times g'$, which defines the multiplication law in the group $G$, has the property $\varphi(s\tau,s'\tau)=(s+s')\tau$ and is uniquely reconstructed from the structure operator, and hence the Lie algebra of the group $G$.

\textbf{Definition 1.} \textit{A deformed gauge translation group $T^{gH}_{X}$, as well as the deformation with which it was obtained, will be called \textbf{canonical} if, for any two points $x,x'\in X$ there exists a smooth parametric curve $x'=x(s)$ such that}

\begin{equation}\label{17}
(s\tau\times s'\tau)_{x}=\varphi_{x}(s\tau_{x},s'\tau_{x'})=(s+s')\tau_{x},
\end{equation}
where $\tau_{x'}=h_{x'}\langle\tilde{\tau}_{x'}\rangle$ and $\tilde{\tau}_{x'}=\dot{x}'$.

Differentiating definition (\ref{17}) with respect to $s$ and $s'$ at zero simultaneously, we obtain the equation of the given curve $\dot{\tau}_{x}+\gamma_{x}\langle\tau_{x},\tau_{x}\rangle=0$, or taking into account (\ref{16}) $\ddot{x}+\Gamma_{x}\langle\dot{x},\dot{x}\rangle=0$. Thus, for an arbitrary torsion-free space with affine connection, the structure of which is given by the action of the group $T^{gH}_{X}$, \textit{such a curve exists and is a geodesic in the affine parameterization} in the space of affine connection defined by the group $T^{gH}_{X}$.
	
\textbf{Criteria for canonicity.} \textit{Along the geodesic connecting two arbitrary points $x,x'\in X$, the auxiliary functions of the canonical group $T^{gH}_{X}$ satisfy the equations
}
\begin{equation}\label{18}
\mu_{x}(s\tau_{x})\langle\tau_{x}\rangle=\tau_{x}+s\gamma_{x}\langle\tau_{x},\tau_{x}\rangle,
\end{equation}
\begin{equation}\label{19}
\lambda_{x}(s\tau_{x})\langle\tau_{x'}\rangle=\tau_{x},
\end{equation}

Equation (\ref{18}) is obtained by differentiating definition (\ref{17}) with respect to $s$, and (\ref{19}), with respect to $s'$ at zero.

Taking into account the equality $\rho_{x}\langle l',l,l\rangle=\partial^{2}\mu_{x}(0)\langle l,l,l'\rangle$ and twice differentiating Eq. (\ref{18}) with respect to $s$ at zero, we obtain $\rho_{x}\langle l,l',l''\rangle+cycl(l,l',l'')=0$, which, using identity (\ref{14}) and definition (\ref{13}), allows us to establish the following relation $\rho_{x}\langle l',l,l\rangle=\frac{2}{3}R_{x}\langle l,l,l'\rangle$, which holds for the canonical groups $T^{gH}_{X}$.

By means of deformation mappings, the auxiliary mapping $\lambda_{x}$, based on the first formula in (\ref{2}), is defined as follows:

\begin{equation}\label{20}
\lambda_{x}(t)=\tilde{\partial}H_{x}(\tilde{t})|_{\tilde{t}=K_{x}(t)}\circ k_{x'},
\end{equation}
as a result of which criterion (\ref{19}) can be represented in the form $\tilde{\partial}H_{x}(\tilde{t}_{x})\langle\tilde{\tau}_{x'}\rangle=\dot{H}_{x}(\tilde{t}_{x})=\tau_{x}$, whence it follows

\textbf{Proposition 1}. \textit{Along the geodesic connecting arbitrary two points separated by a finite interval $\tilde{t}_{x}=x(s)-x(0)$, the corresponding parameter $t_{x}=H_{x}(\tilde{t}_{x})$ of the canonical group $T^{gH}_{X}$ is proportional to the initial vector of the geodesic: $t_{x}=s\tau_{x}$.}

Thus, in the case of a canonical deformation, the function $x'=x+K_{x}(s\tau_{x})$ is a solution to the geodesic equation with initial data $x$, $\dot{x}=k_{x}\langle\tau_{x}\rangle$. The parameters $t_{x}$ of the canonical group $T^{gH}_{X}$ are elements of the tangent bundle $TX$ and are similar to geodesic Riemannian coordinates outgoing from the point $x$.

By virtue of criterion (\ref{19}), the functions

\begin{equation}\label{21}
u(x',\tilde{\tau}_{x'}):=\partial_{x'}H_{x}(x'-x)\langle\tilde{\tau}_{x'}\rangle
\end{equation}
are the first integrals of the geodesic equations $\dot{x}'=\tilde{\tau}_{x'}$, $\dot{\tilde{\tau}}_{x'}=
-\Gamma_{x'}\langle\tilde{\tau}_{x'},\tilde{\tau}_{x'}\rangle$ and, therefore, satisfy the equation

\begin{equation}\label{22}
\partial_{x'}u(x',\tau_{x'})\langle\tilde{\tau}_{x'}\rangle-\partial_{\tau_{x'}}u(x',\tau_{x'})
\langle\Gamma_{x'}\langle\tilde{\tau}_{x'},\tilde{\tau}_{x'}\rangle\rangle=0
\end{equation}
with boundary conditions $u(x,\tilde{\tau}_{x}):=k_{x}\langle\tilde{\tau}_{x}\rangle$, whose characteristics are geodesics. Equation (\ref{22}), taking into account expression (\ref{21}), leads to the following equation for the deformation functions:

\begin{equation}\label{23}
(\partial^{2}_{x'}H_{x}(x'-x)-\partial_{x'}H_{x}(x'-x)\circ
\Gamma_{x'})\langle\tilde{\tau}_{x'},\tilde{\tau}_{x'}\rangle=0
\end{equation}
with boundary conditions

\begin{equation}\label{24}
H_{x}(0)=0,\quad\partial_{x'}H_{x}(0)=h_{x}.
\end{equation}
The Cauchy problem (\ref{23}), (\ref{24}) has a unique solution.

\textbf{Proposition 2.} \textit{From the deformation coefficients $h_{x}$ and the connection $\Gamma_{x}$, the map of the canonical deformation $H_{x}(\tilde{t}_{x})$, and hence the canonical group $T^{gH}_{X}$, are uniquely determined.}

We especially note that the condition of canonicity does not impose restrictions on the mappings $h_{x}$ and $\Gamma_{x}$.

\section{Group of Riemannian translations}

Since all vectors tangent to the geodesic are parallel to each other, equation (\ref{19}) describes the finite parallel transport of the tangent vector $\tau_{x'}$ along the geodesic to the initial point $x$:

\begin{equation}\label{25}
\tau_{x||}=\tau'_{x}:=\lambda_{x}(s\tau_{x})\langle\tau_{x'}\rangle.
\end{equation}
Since $\lambda_{x}(s\tau_{x})\cong1+s\tau_{x}\cdot\gamma_{x}$, this agrees with the parallel transport rule (\ref{15}) of an arbitrary vector with an infinitesimal displacement, and takes place due to the fulfillment of the composition law for transports of vectors tangent to geodesics along the same geodesics: $\lambda_{x}((s+s')\tau_{x})\langle\tau_{x''}\rangle=
\lambda_{x}(s\tau_{x})\langle\lambda_{x'}(s'\tau_{x'})\langle\tau_{x''}\rangle\rangle$, due to which the finite parallel transport of the tangent vector (\ref{25}) coincides with the integral sequence of infinitesimal transports (as classically the finite parallel transport is defined).

Note that for vectors $\theta$ transverse to the tangents to the geodesic, a similar relation $\lambda_{x}((s+s')\tau_{x})\langle\theta_{x''}\rangle=
\lambda_{x}(s\tau_{x})\langle\lambda_{x'}(s'\tau_{x'})\langle\theta_{x''}\rangle\rangle$ in a curved space not fulfilled. Indeed, since the vector $\theta$ are arbitrary, this relation with help of equation (\ref{5}) leads to the equality $\partial'\varphi_{x}(s\tau_{x},s'\tau_{x'})=\lambda_{x}(s\tau_{x})$, which immediately implies $\rho_{x}=0$, and hence $R_{x}=0$. Thus, in a curved space, the finite $\lambda$-\textit{transport} of the vector $\theta$ transversal to the tangent to the geodesic, defined by the formula $\theta'_{x}:=\lambda_{x}(s\tau_{x})\langle\theta_{x'}\rangle$, is not a parallel transport along geodesic, despite the fact that for infinitesimal displacements (\ref{15}) it it is, and this is precisely due to the fact that in this case the compositional law of transport is violated. In this case, the curvature of the Riemannian space acts as a measure of the deviation of geodesics.

According to Riemann, “\textit{measurement consists in the sequential attachment of comparable quantities, therefore the possibility of measurements is due to the presence of some way to transport one quantity, taken as a unit of scale, to another quantity}" \cite{1}. This is exactly what the canonical group $T^{gH}_{X}$ does. Indeed, by virtue of the canonicity, the operation of n-fold transport of the tangent vector $\tau_{x}$ along the geodesic gives

\[
x_{n}=x+K_{x}(\tau_{x})+K_{x_{1}}(\tau_{x_{1}})+\ldots+K_{x_{n-1}}(\tau_{x_{n-1}})=x+K_{x}(n\tau_{x}).
\]
Thus, the geodesic connecting the points $x$ and $x'=x+K_{x}(s\tau_{x})$ contains $s$ vectors $\tau$ when they are transported along themselves, and the affine parameter of the geodesic $s$ gives the measure of its length on the scale of the vector $\tau$.

Specifying the length of a vector $\tau_{x}\in T$, regardless of its direction and location $\tau^{2}_{x}:=\eta\langle\tau_{x},\tau_{x}\rangle$, i.e., specifying the Euclidean metric $\eta$ in the space $T$, allows one to compare the lengths of arbitrary geodesics, which implements the Riemann length principle. Moreover, in this way the structure of a \textit{Riemannian space with metric} $g_{x}\langle\tilde{\tau}_{x},\tilde{\tau}_{x}\rangle:=\eta\langle h_{x}\langle\tilde{\tau}_{x}\rangle,h_{x}\langle\tilde{\tau}_{x}\rangle\rangle$ is setting on $X$. The vector $\tilde{\tau}_{x}$, tangent to the geodesic, acts as a scale unit, and therefore it is natural to require that its length be invariable when transported along the geodesic.

\textbf{Definition 2.} \textit{The canonical gauge translation group, for which the condition of preserving the length of the unit of the scale fulfilled when it is transported along itself $\eta\langle\tau_{x},\tau_{x}\rangle=\eta\langle\tau_{x'},\tau_{x'}\rangle$, is called the \textbf{group of Riemannian translations} and denoted as $RT$.}

Taking into account criterion (\ref{19}), this condition can be represented in the form $G_{x}(s\tau_{x})\langle\tau_{x'},\tau_{x'}\rangle=
\eta\langle\tau_{x'},\tau_{x'}\rangle$, where $G_{x}(t)\langle\tau,\tau\rangle:=\eta\langle \lambda_{x}(t)\langle\tau\rangle,\lambda_{x}(t)\langle\tau\rangle\rangle$, or, using expression (\ref{20}), as

\begin{equation}\label{26}
\eta\langle\partial_{x'}H_{x}(x'-x)\langle\tilde{\tau}_{x'}\rangle,
\partial_{x'}H_{x}(x'-x)\langle\tilde{\tau}_{x'}\rangle\rangle=g_{x'}\langle\tilde\tau_{x'},\tilde\tau_{x'}\rangle.
\end{equation}
Differentiation of equation (\ref{26}) with respect to $x'$ at $x'\rightarrow x$ gives the condition of compatibility the connection with the metric

\begin{equation}\label{27}
2g_{x}\langle\tilde{\tau},\Gamma_{x}\langle\tilde{t},\tilde{\tau}\rangle\rangle=
\partial g_{x}\langle\tilde{t},\tilde{\tau},\tilde{\tau}\rangle,
\end{equation}
or taking into account (\ref{16})

\begin{equation}\label{28}
\eta\langle\tau,\gamma_{x}\langle t,\tau\rangle\rangle=0,
\end{equation}
where all vectors are taken at the point $x$. Condition (\ref{27}), together with the fact that the mapping $\Gamma_{x}$ is symmetric, allows us to express it in terms of the derivative of the metric $\partial g_{x}$:

\[
g_{x}\langle\tilde{t},\Gamma_{x}\langle\tilde{\tau},\tilde{\tau}\rangle\rangle=
\partial g_{x}\langle\tilde{\tau},\langle\tilde{t},\tilde{\tau}\rangle\rangle-
\frac{1}{2}\partial g_{x}\langle\tilde{t},\langle\tilde{\tau},\tilde{\tau}\rangle\rangle,
\]
and condition (\ref{28}) together with (\ref{10}) leads to an expression for the mapping $\gamma_{x}$ of the group $RT$ in terms of its structure operator

\[
\eta\langle\tau,\gamma_{x}\langle t,\tau'\rangle\rangle=
\frac{1}{2}(\eta\langle t,C_{x}\langle\tau,\tau'\rangle\rangle+
\eta\langle\tau,C_{x}\langle t,\tau'\rangle\rangle-\eta\langle\tau',C_{x}\langle t,\tau\rangle\rangle),
\]
whence it follows that in this case the connection $\Gamma_{x}$ is formed by the Christoffel symbols, and $\gamma_{x}$ by the Ricci rotation coefficients and, therefore, are uniquely determined by the frame fields $e_{x}$ which defined on $X$ by the action of the group $RT$.

Condition (\ref{26}) does not impose restrictions on the deformation coefficients $h_{x}$ (and hence the fields $e_{x}$); therefore, taking into account Theorem 1 and Proposition 2, we arrive at the following statement.

\textbf{Theorem 2.} \textit{The group of Riemannian translation $RT$, by its action on the manifold $X$, defines the structure of a Riemannian space and an orthonormal frame field $e_{x}$, by which it is uniquely determined. An arbitrary Riemannian structure on a manifold $X$ can be defined in this way.}

The requirement that the length of the vector $\tau$ tangent to the geodesic be preserved under finite $\lambda$-transport does not provide a similar property for an arbitrary vector $\theta$.

\textbf{Proposition 3.} \textit{The requirement that the length of an arbitrary vector be preserved under a finite $\lambda$-transport $G_{x}(s\tau)=\eta$ is a condition of flat space.}

The statement follows from the relation

\[
\frac{d^{2}}{ds^{2}}G_{x}(s\tau)\langle\theta,\theta\rangle|_{s=0}=
\frac{4}{3}\eta\langle\theta,R_{x}\langle\tau,\tau,\theta\rangle\rangle,
\]
which holds for the group RT \cite{12}.

Thus, in spite of the fact that the group $RT$ contains complete information on the geometrical structure of the Riemannian space and defines it by its action on $X$, it leaves no invariants either in $X$ or in $TX$, which does not allow us to consider group $RT$ as the principal group of Riemannian space. The group of Riemannian translation $RT$ is the source of geometric objects in Riemannian space: its generators form an orthonormal frame field, which defines a metric on the manifold $X$; parameters are vector fields on $X$. The group $RT$ allows one to measure the lengths of arbitrarily located curves by moving a unit scale along them and thus implements the \textit{principle of length} (B. Riemann).

\section{Group of motions of tangent bundle of Riemannian space}

The action of group of Riemannian translations $RT$ defines on the manifold $X$ not only a Riemannian structure, that is, a metric $g_{x}$, but also a fixed field of orthonormal frames $e_{x}=k_{x}\cdot\partial_{x}$. Therefore, the same Riemannian structure on $X$ is given by all groups $RT$ whose frame fields are related by transformations from the gauge rotation group $\tilde{r}_{x}\in R^{g}$:

\begin{equation}\label{29}
e'_{x'}=\tilde{r}^{-1}_{x}\cdot e_{x'},\quad\eta\langle r_{x}\langle\tau\rangle,r_{x}\langle\tau\rangle\rangle=
\eta\langle\tau,\tau\rangle\ \forall x\in X,\ \tau\in T,
\end{equation}
at which the deformation coefficients of the group $RT$ transform according to the law $h'_{x'}=\tilde{r}_{x}\circ h_{x'}$, which leaves the metric $g_{x}$ unchanged. Transformations (\ref{29}) are taken non-local ($x'\neq x$!), which determines a moving frame (according to E. Cartan): here $e'_{x'}$ is interpreted as a frame transported from point $x$ to point $x'$ and arbitrarily turned by the value $\tilde{r}_{x}^{-1}$, given at the point $x$, where the translation value $\tilde{t}_{x}=x'-x$ given, too.

The set of such transformations forms a generalized gauge group $G^{g}_{X}=T^{g}_{X}\times)R^{g}$ \cite{8}, which has the structure of a semidirect product $\times)$ of its subgroups, with parameters $\tilde{\vartheta}_{x}=(\tilde{t}_{x},\tilde{r}_{x})$ and the multiplication law $(\tilde{\vartheta}\times\tilde{\vartheta}')_{x}=:\tilde{\Phi}_{x}(\tilde{\vartheta}_{x},\tilde{\vartheta}'_{x'})$:

\begin{equation}\label{30}
(\tilde{\vartheta}\times\tilde{\vartheta}')^{T}_{x}=\tilde{\Phi}_{x}^T(\tilde{\vartheta}_{x},\tilde{\vartheta}'_{x'}):=\tilde{t}_{x}+\tilde{t}'_{x'},
\end{equation}
\begin{equation}\label{31}
(\tilde{\vartheta}\times\tilde{\vartheta}')^{R}_{x}=\tilde{\Phi}_{x}^R(\tilde{\vartheta}_{x},\tilde{\vartheta}'_{x'}):=\tilde{r}_{x}\circ\tilde{r}'_{x'},
\end{equation}
where $x'=x+\tilde{t}_{x}$ defines the action of the group $G^{g}_{X}$ on $X$, and the indices $T$ and $R$ indicate the canonical projection of the elements of the group $G^{g}_{X}$ onto its factors $T^{g}_{X}$ and $R^{g}$.

The condition $\tau'_{x}\cdot e'_{x'}=\tau_{x'}\cdot e_{x'}$, which follows from the interpretation of transformations (\ref{29}) as frame motions, determines the action of the group $G^{g}_{X}$ on the parameters of the group $RT$:

\begin{equation}\label{32}
\tau'_{x}=\tilde{r}_{x}\langle\tau_{x'}\rangle,
\end{equation}
that is, in the tangent bundle $TX$. Thus, the group $G^{g}_{X}$ describes \textit{nonlocal automorphisms of the group} $RT$ that leave the Riemannian structure it defines on $X$ unchanged.

Transformations (\ref{29}) and (\ref{32}) in the absence of rotations $\tilde{r}_{x}=1$ describe the componentwise identification of vectors given at different points, which corresponds to the case of flat space and does not agree with the Riemannian structure defined on $X$ by the action of the group $RT$. To implement such a match, we deform the group $G^{g}_{X}$.

Let the \textit{deformed group} $G^{g\bar{H}}_{X}=T^{g\bar{H}}_{X}\times)R^{g}$ be parameterized by the pairs $\vartheta_{x}=(t_{x},r_{x})$ in such a way that

\begin{equation}\label{33}
\tilde{t}_{x}=\bar{K}_{x}(t_{x}),
\end{equation}
\begin{equation}\label{34}
\tilde{r}_{x}=r_{x}\circ\pi_{x}(t_{x}),
\end{equation}
where $\bar{K}_{x}=\bar{H}^{-1}_{x}$ and $\pi_{x}(t)\in R^{g}$, that is,

\begin{equation}\label{35}
\eta\langle \pi_{x}(t)\langle\tau\rangle,\pi_{x}(t)\langle\tau\rangle\rangle=
\eta\langle\tau,\tau\rangle,
\quad\pi_{x}(0)=1\ \forall x\in X,\ t,\tau\in T.
\end{equation}

Formula (\ref{33}) describes an independent deformation of the subgroup $T^{g}_{X}\subset G^{g}_{X}$ to $T^{g\bar{H}}_{X}$ using the deformation map $\bar{H}_{x}$, and (\ref{34}), being substituted in (\ref{32}), ensures the rotation $\pi_{x}(t_{x})$ vectors under translations:

\begin{equation}\label{36}
\tau'_{x}=r_{x}\langle \pi_{x}(t_{x})\langle\tau_{x'}\rangle\rangle,
\end{equation}
where $x'=x+\bar{K}_{x}(t_{x})$ defines the action of the group $G^{g\bar{H}}_{X}$ on $X$, and (\ref{36}) on $TX$.

The multiplication law $(\vartheta\times\vartheta')_{x}=:\Phi_{x}(\vartheta_{x},\vartheta'_{x'})$ of the group $G^{g\bar{H}}_{X}$ is determined by its isomorphism (\ref{33}), (\ref{34}) to the group $G^{g}_{X}$ and the multiplication law (\ref{30}), (\ref{31}) of the group $G^{g}_{X}$:

\[
(\vartheta\times\vartheta')^{T}_{x}=\Phi_{x}^T(\vartheta_{x},\vartheta'_{x'}):=\bar{\varphi}_{x}(t_{x},t'_{x'})=\bar{H}_{x}(\bar{K}_{x}(t_{x})+
                                        \bar{K}_{x'}(t'_{x'})),
\]
\begin{equation}\label{37}
(\vartheta\times\vartheta')^{R}_{x}=\Phi_{x}^R(\vartheta_{x},\vartheta'_{x'}):=r_{x}\circ\pi_{x}(t_{x})\circ r'_{x'}\circ\pi_{x'}(t'_{x'})\circ\pi^{-1}_{x}(\bar{\varphi}_{x}(t_{x},t'_{x'})) .
\end{equation}

Transformation (\ref{36}) in the case of pure displacement, i.e., with $\vartheta_{x}=(t_{x},1)$, is called the $\pi$-\textit{transport}. In the infinitesimal case, it takes the form $\tau'_{x}=\tau_{x}+t_{x}\cdot\bar{\nabla}_{x}\tau_{x}$, where $\bar{\nabla}_{x}:=\bar{e}_{x}+\bar{\gamma}_{x}$, $\bar{e}_{x}:=\bar{k}_{x}\cdot\partial_{x}$, $\bar{k}_{x}:=\partial\bar{K}_{x}(0)$, $\bar{\gamma}_{x}:=\partial\pi_{x}(0)$, and corresponds to the parallel transport of the vector $\tau_{x}$ from the point $x'=x+\bar{k}_{x}\langle t_{x}\rangle$ to the point $x$ in the affinely connected space with the connection $\bar{\gamma}_{x}$ given in the frame $\bar{e}_{x}$. In this sense, this geometric structure is defined on $X$ by the action (\ref{36}) of the group $G^{g\bar{H}}_{X}$ in $TX$.

We require that the \textit{consistency condition} $\bar{\nabla}_{x}=\nabla_{x}$ (or $\bar{e}_{x}=e_{x}$, $\bar{\gamma}_{x}=\gamma_{x}$) be satisfied, which ensures coincidence at the infinitesimal level of $\pi$-transports with $\lambda$-transports, and hence the coincidence of the geometric structure given by the group $G^{g\bar{H}}_{X}$, to the structure given by the original group $RT$, whose nonlocal automorphisms are described by the group $G^{g\bar{H}}_{X}$. Note that, due to requirement (\ref{35}), condition of metric compatibility (\ref{28}) is satisfied automatically, and the condition $\bar{\gamma}_{x}=\gamma_{x}$, due to the accepted equality $\bar{e}_{x}=e_{x}$, is reduced to the nontorsionity condition

\[
C_{x}\langle l,l'\rangle:=\bar{\gamma}_{x}\langle l,l'\rangle-(l\leftrightarrow l')
\]

\textbf{Definition 3.} \textit{The deformed group $G^{g\bar{H}}_{X}$ of nonlocal automorphisms of the group of Riemannian translations $RT$ preserving the Riemannian structure defined on $X$ by the action of the group $RT$, for which the \textbf{consistency condition} $\bar{\nabla}_{x}=\nabla_{x}$ is satisfied, will be called the \textbf{group of motions of tangent bundle of Riemannian space} and denoted as $DR$.}

In terms of the mappings $L_{x}(\vartheta):=r\circ\pi_{x}(t)$, the compositional law of transformations (\ref{36}) $\tau'_{x}=L_{x}(\vartheta)\langle\tau_{x'}\rangle$, written for constant parameters $\vartheta=(t,r)$, takes the form

\begin{equation}\label{38}
L_{x}(\Phi_{x}(\vartheta,\vartheta'))=L_{x}(\vartheta)\circ L_{x'}(\vartheta'),
\end{equation}
where $x'=x+\bar{K}_{x}(t)$. Introducing the notation

\[
\Lambda_{x}(\vartheta):=\partial_{\vartheta'}\Phi_{x}(\vartheta,\vartheta')|_{\vartheta'=(0,1)},\quad E_{x}:=\partial_{\vartheta}L_{x}(\vartheta)|_{\vartheta=(0,1)},
\]
where $\partial_{\vartheta}:=\partial/\partial\vartheta$, and differentiating equation (\ref{38}) with respect to $\vartheta'$ at $\vartheta'=(0,1)$, we obtain an analogue of the Lie equation for the group $DR$ of transformations (\ref{36}):

\begin{equation}\label{39}
\Lambda_{x}(\vartheta)\cdot\partial_{\vartheta}L_{x}(\vartheta)=L_{x}(\vartheta)\circ E_{x'}.
\end{equation}

The integrability condition for equation (\ref{39}) is the equation

\begin{equation}\label{40}
[D_{x},D_{x}]=\Sigma_{x}\cdot D_{x},
\end{equation}
written by us in terms of the commutator of the generators $D_{x}$ of the group $DR$ of transformations (\ref{36}), which are determined by the relation

\[
D_{x}\tau_{x}:=\partial_{\vartheta}(L_{x}(\vartheta_{x})\langle\tau_{x'}\rangle)|_{\vartheta=(0,1)},
\]
and its structure operator

\begin{equation}\label{41}
\Sigma_{x}\langle\xi,\xi'\rangle:=\partial_{\vartheta}\Lambda_{x}(\vartheta)\langle\xi,\xi'\rangle|_{\vartheta=(0,1)}-
(\xi\leftrightarrow\xi'),
\end{equation}
where in this case $\xi$ and $\xi'$ are the vectors of the space tangent to the unit of the group $G=T\otimes R$. Equation (\ref{40}) generalizes the \textit{Maurer – Cartan equation} of the theory of finite-dimensional Lie groups of transformations for the group $DR$ of transformations (\ref{36}).

The generators $D_{x}$ split into translation generators, which coincide with the covariant derivatives $\nabla_{x}$ of vector fields, and also $x$-independent generators $A$ of rotations of vectors in the Euclidean space $T$.

A special case of (\ref{40}) is the equation

\begin{equation}\label{42}
[\nabla_{x},\nabla_{x}]=\Sigma_{x}{}^{T}_{TT}\cdot\nabla_{x}+\Sigma_{x}{}^{R}_{TT}\cdot A,
\end{equation}
for the commutator of translation generators, where the subscript $T$ means the restriction of the mapping to the subspace $T$. Direct calculations using the defining formula (\ref{41}) give the following expressions for the components of the structure operator of the group $DR$: $\Sigma_{x}{}^{T}_{TT}=C_{x}$, $\Sigma_{x}{}^{R}_{TT}=R_{x}$, as a result of which equation (\ref{42}) takes the form of the structural equation of a Riemannian space

\begin{equation}\label{43}
[\nabla_{x},\nabla_{x}]=C_{x}\cdot\nabla_{x}+R_{x}
\end{equation}
with the curvature tensor $R_{x}$ in the frame $e_{x}$, whose anholonomic object is $C_{x}$ (here it is taken into account that $R_{x}=R_{x}\cdot A$). Equation (\ref{43}) using equation (\ref{9}) $[e_{x},e_{x}]=C_{x}\cdot e_{x}$ is reduced to (\ref{12}).

\textit{Jacobi identity} for the group $DR$

\[
D_{x}\Sigma_{x}\langle\xi,\xi',\xi''\rangle+
\Sigma_{x}\langle\xi,\Sigma_{x}\langle\xi',\xi''\rangle\rangle+cycl(\xi,\xi',\xi'')=0
\]
is the integrability condition for the Maurer – Cartan equation (\ref{40}), and when it is restricted to the translation $\xi\rightarrow l$ reduces to the \textit{Bianchi identity}.

Thus, it is proved

\textbf{Theorem 3.} \textit{The translation generators of the action of the group of motions of the tangent bundle of the Riemannian space $DR$ in its tangent bundle $TX$ are covariant derivatives of vector fields $\nabla_{x}=e_{x}+\gamma_{x}$, the structure operator of the group $DR$ has as components the object of anholonomity $C_{x}$ of the orthonormal frame field $e_{x}$ and the curvature tensor $R_{x}$ in it, and the structure equation of the Riemannian space $[\nabla_{x},\nabla_{x}]=C_{x}\cdot\nabla_{x}+R_{x}$ is a necessary condition for the existence of the group $DR$ which define by its action in the bundle $TX$ this structure on $X$.}

We see that both $DR$ and $RT$ define the structure of a Riemannian space on $X$, but in two different ways. If for group $DR$ an infinitesimal action in the tangent bundle $TX$ is sufficient and structure equation (\ref{43}) appears as a component of the Maurer – Cartan equation which ensure the existence of the group $DR$, then for group $RT$ it is required to specify its action in $X$ already with an accuracy of at least the second order in displacements, which generates both the bundle $TX$ itself and the action of the group $RT$ in it, and the structure equation of the Riemannian space (\ref{12}) appears here already after the differentiation of the Maurer – Cartan equation for the group $RT$. The reason for this is that the groups $DR$ and $RT$ describe different aspects of the curvature of the Riemannian space, namely, in the case of group $DR$, curvature acts as a measure of the noncommutativity of its generators (or the measure of the rotation of the vector when it is swept along a closed contour), and in the case of group $RT$, curvature acts as a measure of geodesic deviation.

We emphasize that in the proposed group-theoretic approach, the curvature is not added to the structural equation of a flat space for reasons of anholonomicity of curved space, as do E. Cartan, but is a characteristic of a group that defines a given geometric structure with variable curvature. In the case of the group $DR$, this characteristic is its structure operator, which is determined by the antisymmetric part of the coefficients in the second order of the expansion of the multiplication law of the group $DR$ with respect to parameters, and in the case of group $RT$, it is the curvature operator, which is determined by the antisymmetric part of the coefficients already in the third order of the expansion of the multiplication law of the group $RT$. Thus, the anholonomicity of the space appears upon deformation of the generalized gauge group of the corresponding flat space, as an isomorphic transition to the deformed group, which \textit{preserves the group structure of the transformations used, including their associativity}.

In the orthonormal frame bundle (in principal bundle $RX$), the group of motions of the tangent bundle of the Riemannian space $DR$ acts by the formula

\[
e'_{x'}=\pi^{-1}_{x}(t_{x})\cdot r^{-1}_{x}\cdot e_{x'},
\]
This formula follows from (\ref{29}) (taking into account deformation (\ref{34})), determines the movable frame of E. Cartan \cite{10} and allows one to define connection by means of a movable frame (with r = 1):

\[
\nabla_{e_{x}}e_{x}:=\lim_{t\rightarrow0}\frac{e_{x}-e'_{x}}{t}=\gamma_{x}\cdot e_{x}.
\]

Since $r_{x},\pi_{x}\in R^{g}$, transformations from the group of motions of the tangent bundle of the Riemannian space $DR$ \textit{preserve the lengths of vectors and the angles between them}, and therefore allow one to establish the equality of geometric figures by superimposing them, even if these figures are located at different points of the space $X$. Here, under geometric figures we mean configurations of vectors defined in tangent spaces. Thus, $DR$ realizes the \textit{principle of equality} and \textit{is the principal group of Riemannian space} (according to F. Klein).
	
\section{Group of parallel transports of Riemannian space}

Since the infinitesimal action of the group $DR$ is sufficient to define the geometric structure, the same Riemannian structure is given by the set of groups $DR$ differing in higher orders of action in the bundles $TX$ or $RX$. On the other hand, the question arises - what does the finite $\pi$-transport mean? Along what curve is it carried out? These questions are similar to the questions about the meaning of finite $\lambda$-transports in the case of the deformed group of diffeomorphisms $T^{gH}_{M}$, which were successfully clarified for it using the concept of canonicity.

For the group $DR$, which is broader than the group $T^{gH}_{M}$, the concept of canonicity is somewhat broadened, but equally fruitful.

\textbf{Definition 4.} \textit{The group of motions of the tangent bundle of the Riemannian space $DR$ will be called the \textbf{group of parallel transports of Riemannian space} and denoted as $DP$ if it satisfies the extended \textbf{canonicity condition}, namely, for any two points $x,x'\in X$, there exists a smooth parametric curve $x'=x(s)$, $x=x(0)$, such that}

\[
((s\tau,1)\times(s'\tau,1))_{x}=(\bar{\varphi}_{x}(s\tau_{x},s'\tau_{x'}),1)=((s+s')\tau_{x},1),
\]
\textit{where $\tau_{x'}=h_{x'}\langle\tilde{\tau}_{x'}\rangle$, $\tilde{\tau}_{x'}=\dot{x}'$.}
	
This definition, firstly, requires the subgroup $T^{g\bar{H}}_{M}\subset DP$ to be canonic, and, due to the consistency condition $\bar{e}_{x}=e_{x}$ and Theorem 2, the coincidence of the deformation mappings $\bar{H}_{x}=H_{x}$, and hence the coincidence of the entire group $T^{g\bar{H}}_{M}$ with the original group of Riemannian translation $RT$, which also implies the existence of the desired curves and their coincidence with the geodesics given by the group $RT$. Second, as follows from formula (\ref{37}), it leads to the fulfillment of the composition law for $\pi$-transports of arbitrary vectors along geodesics: $\pi_{x}((s+s')\tau_{x})=\pi_{x}(s\tau_{x})\circ\pi_{x'}(s'\tau_{x'})$. Since for infinitesimal displacements, the $\pi$-transport is a parallel transport of vectors, the fulfillment of this law leads to the fact that the finite $\pi$-transport of arbitrary vectors coincides with the integral sequence of their infinitesimal parallel transports along the geodesic connecting points $x$ and $x'$  (which determined the name of the group $DP$ ).

Due to the consistency condition, in the group of parallel transports $DP$, as in all groups $DR$, the infinitesimal $\pi$-transports and $\lambda$-transports of arbitrary vectors coincide, while the finite transports in a curved space are different. But only in the group $DP$ for vectors tangent to geodesics, both finite $\pi$-transports and $\lambda$-transports lead to the same results.

To find the equation defining the mapping $\pi_{x}$, note that the functions

\begin{equation}\label{44}
u(x',\theta_{x'}):=\pi_{x}(H_{x}(x'-x))\langle\theta_{x'}\rangle,
\end{equation}
due to the fact that along geodesics connecting the pairs of points $x$ and $x'$ $u(x',\theta_{x'})=\theta_{x}=const$, are the first integrals of the system of equations $\dot{x}'=k_{x}\langle\tau_{x'}\rangle$, $\dot{\theta}_{x'}=-\gamma_{x'}\langle\tau_{x'},\theta_{x'}\rangle$, and therefore satisfy the partial differential equation:

\[
k_{x}\cdot\partial_{x'}u(x',\theta_{x'})\langle\tau_{x'}\rangle-
\partial_{\theta_{x'}}u(x',\theta_{x'})\langle\gamma_{x'}\langle\tau_{x'},\theta_{x'}\rangle\rangle=0,
\]
substitution into which expression (\ref{44}) leads to the equation

\[
k_{x'}\cdot\partial_{x'}\pi_{x}(H_{x}(x-x'))-\pi_{x}(H_{x}(x-x'))\circ\gamma_{x'})\langle\tau_{x'},\theta_{x'}\rangle=0,
\]
giving, under the condition $\pi_{x}(0)=1$ and arbitrary initial values $\tau_{x},\theta_{x}\in T$, a unique solution for the map $\pi_{x}$. Thus, the deformation (\ref{33}), (\ref{34}) with the help of which the group of parallel transports $DP$ is constructed, is uniquely determined.

\textbf{Proposition 4.} \textit{For a given orthonormal frame field $e_{x}$, the group of parallel transports of Riemannian space $DP$ is uniquely determined.}

Thus, we can say that the group of parallel transports $DP$ is univocal extension of the group of Riemannian translations $RT$ to the possibility of taking into account local rotations from $R^g$. Uniqueness is achieved, first of all, due to the canonicity condition, which, from a geometric point of view, is reduced to fixing geodesics as curves along which finite transports are carried out, as well as imposing a consistency condition. All groups $DP$ correspond to the same Riemannian structure, if for them frames $e_{x}$ are related by transformations from $R^g$.

\textbf{Theorem 4.} \textit{The group of parallel transports of Riemannian space $DP$, acting on the manifold $X$ as the group of Riemannian translations $RT$ (with the inefficiency kernel $R^g$), defines the structure of a Riemannian space on $X$, making it possible to measure the lengths of arbitrary curves by moving a unit scale, which implements the \textbf{principle of length}.}

\textit{Acting in the tangent bundle $TX$, group $DP$ preserves the lengths of vectors and the angles between them, making it possible to establish the congruence of geometric shapes by superimposing them, which implements the \textbf{principle of equality}.}

So the group of parallel transports of Riemannian space $DP$ combines both principles of the foundations of geometry, thus \textit{overcoming the Riemann-Klein antagonism.}

\section*{Conclusions}

The main results of this article are as follows:

 the group of Riemannian translations $RT$ is the source of geometric objects of Riemannian space and implements the length principle;

 the group of motions of tangent bundle of Riemannian space $DR$, as the group of nonlocal automorphisms of the group $RT$ compatible with $RT$, is the principal group of Riemannian space, which realizes F. Klein's Erlangen program for it - the equality principle;

 the group of parallel transports $DP$, as a canonical group $DR$, containing the group of Riemannian translations $RT$ as a subgroup, combines both principles and coordinates them.


\section*{References}
\begin{thebibliography}{<99>}

\bibitem{1} Riemann, B. \"{U}ber die Hypothesen, Welche der Geometre zu Gruhde Liegen, in: {\it On Foundations of Geometry}, pp. 309--325. Gostekhteorizdat, Moscow (1956) (Russian translation)
\bibitem{2} Klein F. Vergleichende Betrachtungen Uber Neuere Geometrische Forschungen (Erlangen Program), in: {\it On Foundations of Geometry}, pp. 399--434. Gostekhteorizdat, Moscow (1956) (Russian translation)
\bibitem{3} Cartan E. Group Theory and Geometry, in: {\it On Foundations of Geometry}, pp. 438--507. Gostekhteorizdat, Moscow (1956) (Russian translation)
\bibitem{4} Sulanke R. and Wintgen P. Differentialgeometrie und Faserbundel, Veb Deutscher Verlag der Wissenschaften, Berlin (1972)
\bibitem{5} Sabinin L. Methods of Nonassociative Algebra in Differential Geometry, in: Koboyashi S. and Nomizu K. (eds.) {\it Foundations of Differential Geometry} {\bf1}, pp. 293--334. Nauka, Moscow (1981) (in Russian)
\bibitem{6} Samokhvalov S. Group-theoretical Description of Gauge Fields, {\it Theor. Math. Phys.} {\bf76} (1988) 709--717
\bibitem{7} Samokhvalov S. On Specification of Connections in Bundles by the Action of Infinite Lie Groups, {\it Ukrainian Math. J.} {\bf43} (1991) 1599--1603
\bibitem{8} Samokhvalov S. and Reznyk K. Cartan Equation as the Condition of the Existence of an Infinite Group, {\it Bulg. J. Phys.} {\bf76} No.S2 (2006) 309–314
\bibitem{9} Samokhvalov S. Group-theoretic Description of Riemannian Spaces, {\it Ukr. Math. J.} {\bf55} (2003) 1238--1248, arXiv:0704.2967 [math.DG]
\bibitem{10} Cartan E. Riemannian Geometry in Ortonormal Frame, Univ. Press, Moscow (1960) (Russian translation)
\bibitem{11} Ovsyanikov L. Group Analysis of Differential Equations, Nauka, Moscow (1978). (in Russian)
\bibitem{12} Samokhvalov S. Fundamental Group of the Einstein Space, {\it Mat. Modelyuvannya} No.2 (2008) 15–19 (in Ukrainian)

\end{thebibliography}
\end{document}